\documentclass[12pt]{amsart}
\usepackage{amsthm,vmargin}
\usepackage{amsmath,xypic}
\usepackage{hyperref}
\usepackage{epigraph}
\usepackage{graphicx}

\newtheorem{theorem}[equation]{Theorem}      
\newtheorem{corollary}[equation]{Corollary}  
\newtheorem{proposition}[equation]{Proposition}

\theoremstyle{definition}

\theoremstyle{definition}

\theoremstyle{remark}

\theoremstyle{definition}
\newtheorem{remark}[equation]{Remark}

\numberwithin{equation}{section} 

\let\isom=\simeq

\newcommand{\bl}{{\rm Bl}}

\newcommand{\M}{\mathcal{M}}
\newcommand{\Mb}{{\overline{\M}}}
\newcommand{\mydot}{{\small{\bullet}}}

\renewcommand{\int}{\operatorname{int}}
\renewcommand{\O}{{\mathcal O}}
\renewcommand{\P}{{\mathbb P}}

\newcommand{\sG}{{\mathcal G}}
\newcommand{\sS}{{\mathcal S}}
\let\sg=\sG
\begin{document}

\title[]{Ordinarity of configuration spaces and of wonderful compactifications}%
\author{Kirti Joshi}%
\address{Math. department, University of Arizona, 617 N Santa Rita, Tucson
85721-0089, USA.} \email{kirti@math.arizona.edu}
\date{Preliminary Version: \jobname--0.1}

\thanks{}%


\begin{abstract}
We prove the following: (1) if $X$ is ordinary, the
Fulton-MacPherson configuration space $X[n]$ is ordinary for all
$n$; (2) the moduli of stable $n$-pointed curves of genus zero is
ordinary. (3) More generally we show that a wonderful compactification
 $X_\sg$ is ordinary if and only if $(X,\sg)$ is an ordinary building
 set. This implies the ordinarity of many other well-known configuration
 spaces (under suitable assumptions).
\end{abstract}
\maketitle
\epigraph{O Marvelous! what new configuration
 will come next?\\
I am bewildered with multiplicity.}{William Carlos Williams}


\section{Introduction}
In the past few years a number of configuration spaces have been
studied (see \cite{fulton94},\cite{deconcini95}, \cite{ulyanov02},
\cite{hu03},\cite{li06}, \cite{chen06},\cite{kim08}).  This class of
schemes also include the moduli of $n$-pointed stable curves of
genus zero, denoted here by $\Mb_{0,n}$ (for $n\geq 3$). All these
configuration spaces typically arise from an \emph{initial datum},
which usually consists of a collection of closed non-singular
subschemes of a non-singular, projective variety with certain
additional properties--like transversal intersection; as well
combinatorial data such as an integer or a graph. Given an initial
datum, the configuration space associated to it is typically
constructed as a sequence of blowups using the subschemes provided
in the initial datum. Many configuration schemes constructed in the
above references can also be constructed as \emph{wonderful
compactifications} of suitable open varieties constructed from the
initial datum (see \cite{deconcini95,li06}).

Now suppose that $k$ is an algebraically closed field of
characteristic $p>0$. A smooth, projective variety $X/k$ is said to
be \emph{ordinary} if $H^i(X,BW\Omega^j_X)=0$ for all $i,j$. Here
$H^i(X,BW\Omega^j_X)$ are the groups defined in \cite{illusie83}
using the de Rham-Witt complex. The vanishing of these is equivalent
to the vanishing of the Zariski cohomology groups
$H^i(X,B\Omega^j_X)$ for all $i,j$ where for any $j\geq 0$,
$B\Omega^j_X=d\Omega^{j-1}_X$ is the sheaf of locally exact
$j$-forms (see \cite{illusie83}). Ordinarity of a variety is a
difficult condition to check in practice as it requires an
understanding of crystalline Frobenius. Here are some examples of
ordinary varieties: projective spaces, Grassmanians, more generally
homogenous spaces $G/P$ for $G$ semisimple, $P$ parabolic subgroup
of $G$; for abelian varieties ordinarity in the  above sense is
equivalent to ordinarity in the usual sense (invertibility of the
Hasse-Witt matrix); that a general abelian variety with a suitable
polarization is ordinary is a nontrivial result of Peter Norman and
Frans Oort \cite{norman80}; that a general complete intersection in
projective space is ordinary is a delicate result of Luc Illusie
(see \cite{illusie90}).

Our remark in this note is that \emph{a configuration scheme (of the
above type), or more generally a wonderful compactification, arising
from an initial datum is ordinary if and only if it arises from an
ordinary initial datum} (see Theorem~\ref{the:main} and
Corollary~\ref{cor:main}). In particular we prove that the following
schemes are ordinary: (1) if $X$ is a smooth, ordinary, projective
variety and let $X[n]$ be the configuration space of
Fulton-MacPherson (see \cite{fulton94}) and its generalizations (see
\cite{kim08}). The scheme $X[n]$ is a compactification of stable
configurations of $n$-points of $X$. (2) $\Mb_{0,n}$,  the moduli
space of $n$-pointed stable curves of genus zero (\cite{keel92}).
(3) The compactification $X\langle n\rangle$ of Ulyanov
\cite{ulyanov02}, (4) the compactification of Kuiperberg-Thurston,
\cite{li06}, (5) the spaces $T_{d,n}$ of stable, pointed, rooted
trees of \cite{chen06}, (6) the compactification of  open varieties
due to Yi Hu (see \cite{hu03}).

The proof is not difficult but as all of these configuration schemes
play an important role in many areas of algebraic geometry, so their
properties in positive characteristic are not without interest, and
hence worth recording.

This note grew out of our attempt to answer a question raised by
Indranil Biswas (unfortunately we cannot answer his question--see
Remark~\ref{biswas-rem} for more on this). It is a pleasure to thank
him for many conversations about his question. We thank Ana-Maria
Castravet for many conversations about $\Mb_{0,n}$, and especially
pointing out the constructions of \cite{keel92,kapranov93}.

\section{Preliminaries}
Let $k$ be a perfect field of characteristic $p>0$. Let
$W\Omega^\mydot_X$ be the de Rham-Witt complex of $X$. Let
$H^i(X,W\Omega_X^j)$ (for $i+j\leq \dim(X)$) be the de Rham-Witt
cohomology groups. We say that $X$ is ordinary if
$H^i(X,BW\Omega^j_X)=0$ for $i,j\geq 0$ (as a convention we declare the empty scheme to be ordinary).
This is equivalent to the
vanishing of $H^i(X,B\Omega^j_X)=0$ for $i,j\geq0$ where
$B\Omega^j=d(\Omega^{j-1}_X)$ is the sheaf of locally exact
differentials. As we are in characteristic $p>0$, and as $X$ is
smooth of finite type, the sheaf $B\Omega^j_X$ carries a natural
structure of an $\O_X$-module. The condition of ordinarity is
equivalent to the condition:
\begin{equation}
F:H^i(X,W\Omega_X^j)\to H^i(X,W\Omega^j_X)
\end{equation}
is an isomorphism of $W$-modules for all $i,j\geq 0$.  We will use
the following standard results.

For a smooth, projective variety $X$ and $Z\subset X$ a smooth,
closed subscheme, let $\bl_Z(X)$ be the blowup of $X$ along $Z$. For
$Y\subset X$ we write $\tilde{Y}\subset\bl_Z(X)$ for the
\emph{dominant transform} of $Y$ in $\bl_Z(X)$, defined as
$\tilde{Y}=\pi^{-1}(Y)$ if $Y\subset Z$ and  the strict transform of
$Y$ in $\bl_Z(X)$ otherwise.

\begin{proposition}[Ekedahl~\cite{ekedahl85}]\label{ekedahl}
Let $X,Y$ be smooth, projective varieties over $k$. Then $X\times_k
Y$ is ordinary if and only if one of $X,Y$ is ordinary and the other
is Hodge-Witt.
\end{proposition}
\begin{proposition}[Illusie~\cite{illusie90}]\label{illusie}
Let $X$ be a smooth, projective variety over a perfect field $k$.
Let $V$ be a vector bundle on $X$. Let $\P(V)\to X$ be the
associated projective bundle. Then $X$ is ordinary if and only if
$\P(V)$ is ordinary.
\end{proposition}

We need the following version of \cite[Proposition~1.6]{illusie90}:

\begin{proposition}\label{buildingsetlemma}
Let $X$ be a smooth, projective scheme over an algebraically closed
field $k$. Let $Z\subset X$ be a subscheme of $X$ and let $Y\subset
X$ be a smooth, closed subscheme of $X$. Let
$\tilde{Y}\subset\bl_Z(X)$ be the dominant transform of $Y$ in
$\bl_Z(X)$. Then $\tilde{Y}$ is ordinary if and only if $Y,Y\cap Z$
are ordinary.
\end{proposition}
\begin{proof}
We write $\pi:\bl_Z(X)\to X$ for the blowup morphism. Then by
\cite[Proposition~1.6]{illusie90},  $\bl_Z(X)$ is ordinary if and
only if $X,Z$ are ordinary. Next observe that the dominant transform
$\tilde{Z}$ of $Z$ is the exceptional divisor and by \cite[Theorem~8.24(b), page 186]{hartshorne-algebraic}, $\tilde{Z}\to Z$ is
a projective bundle and so by Proposition~\ref{illusie} is
$\tilde{Z}$ is ordinary if and only if $Z$ is ordinary.

Now to prove the assertion. Let $Z\subset X$ be a smooth, proper
subscheme of a smooth, proper $X$. Let $Y\subset X$ be a smooth,
proper subscheme. Let $\tilde{Y}\subset\bl_Z(X)$ be the dominant
transform of $Y$ in $\bl_Z(X)$. We consider several subcases. If $Y$
is a subset of $Z$, then the dominant transform  $\tilde{Y}\to Y$ is
a projective bundle over $Y$ and hence $\tilde{Y}$ is ordinary if
and only if $Y$ is ordinary (by \ref{illusie}). If $Y=Z$, then the
dominant transform $\tilde{Y}$ is the exceptional divisor
$E\subset\bl_Z(X)$. Since $E$ is a projective bundle over $Z$, we
see that $\tilde{Y}=E$ is ordinary if and only if $Z$ is ordinary.
If $Y\not\subset Z$ then we proceed as follows. If $Y\cap
Z=\emptyset$ then $\tilde{Y}\isom Y$ and hence is ordinary as $Y$ is
ordinary. If $Y\cap Z$ is non-empty and by previous considerations,
we may assume that $Y\neq Z$. In this case $\tilde{Y}$ is the blowup
of $Y$ along $Y\cap Z$  and so the result follows from
\cite[Proposition~1.6]{illusie90}. This proves the claim.
\end{proof}

\section{Building sets and wonderful compactification}
Let $X$ be a smooth, projective scheme over an algebraically closed
field $k$. Let $\sS$ be a finite collection of closed, smooth
subschemes of $X$. We say that $\sS$ is an \emph{arrangement} if the
scheme theoretic intersection of any elements of $\sS$ is either
empty or an element of $\sS$.

Let $\sS$ be an arrangement of subschemes of $X$. We say that
$\sG\subset\sS$ is a \emph{building set} if for all
$S\in\sS\backslash\sG$, the minimal elements in $\{G\in\sG:G\supset
S\}$ intersect transversally and their intersection is $S$.

A set of subschemes $\sG$ of $X$ is called a building set if the
collection of all possible intersections of elements of $\sG$ is an
arrangement of subschemes of $X$ and $\sG$ is a building set of this
arrangement.

Let $X$ be a smooth, projective scheme over $k$ and let $\sG$ be a
building set of $X$. Let $X_\sG \subset \prod_{G\in\sG}\bl_G(X)$ be
the closure of $X^o=X\backslash\cup_{G\in\sG}G$. Then we have

\begin{theorem}[\cite{li06}]
Let $X$ be a smooth, projective variety over an algebraically closed
field $k$. Let $\sG$ be a building set of $X$. Then $X_{\sG}$ is a
smooth, projective variety over $k$.
\end{theorem}

The scheme $X_\sG$ is called the wonderful compactification of
$(X,\sG)$.

We say that a building set $\sG$ of $X$ is an \emph{ordinary
building set} if $X$ is ordinary and all the scheme theoretic
intersections of any members of $\sG$ are ordinary (recall that by
our convention empty intersections are also ordinary). We say that
an \emph{arrangement $\sS$ of $X$ is ordinary} if $\sS$ arises from
an ordinary building set.

\begin{theorem}\label{the:main}
Let $X/k$ be a smooth, projective scheme over a perfect field of
characteristic $p>0$. Let $\sg$ be a building set associated to $X$.
Then the wonderful compactification $X_{\sg}$ associated to $X$ is
ordinary if and only if $\sG$ is an ordinary building set.
\end{theorem}
\begin{corollary}\label{cor:main}
Let $X$ be an smooth, projective variety over $k$. Assume that $X$
is ordinary. Then the following schemes associated to $X$ are all
ordinary:
\begin{enumerate}
\item the scheme $X[n]$ of Fulton-MacPherson (see \cite{fulton94})
\item the scheme $X\langle n\rangle$ of Ulyanov (see \cite{ulyanov02})
\item the scheme $X^\Gamma$ of Kuiperberg-Thurston (see \cite{li06})
\item the generalized Fulton-Macpherson configuration scheme $X^n_D,X_D[n]$
(we assume $D$ is a smooth, ordinary subscheme of $X$) of \cite{kim08},
\item the moduli, $\Mb_{0,n}$ (for $n\geq 3$), of
$n$-pointed stable curves of genus zero is ordinary.
\item the scheme of $T^{d,n}$ of stable, $n$-pointed, rooted trees of
$d$-dimensional projective spaces of \cite{chen06}.
\end{enumerate}
\end{corollary}

\begin{proof}[Proof of Theorem~\ref{the:main}]
We recall the details of the construction of $X_\sG$. The
construction is inductively carried out as follows. Let $\sS$ be an
arrangement of $X$ and $\sG$ be a building set of $\sS$. Then assume
that $\sG=\{G_1,\ldots,G_N\}$ is indexed so that $G_i\subset G_j$ if
$i\leq j$. We define $(X_k,\sS^{(k)},\sG^{(k)})$ as follows. For
$k=0$, set $X_0=X,\sS^{(0)}=\sS$, $\sG^{(0)}=\sG$, $G_i^{(0)}=G_i$
for $1\leq i\leq N$. Then $(X_0,\sS^{(0)},\sG^{(0)})$ is ordinary.
Assume by induction that $(X_{k-1},\sS^{(k-1)},\sG^{(k-1)})$ has
been constructed so that $X_{k-1}$ is ordinary and $\sG^{(k-1)}$ is
an ordinary building set for $X_{k-1}$. Then $\sS^{(k-1)}$ consists
of ordinary subvarieties of $X_{k-1}$. Define
$X_k=\bl_{G_k^{(k-1)}}(X_{k-1})$. Then by
Proposition~\ref{buildingsetlemma}, $X_k$ is ordinary if and only if
$X_{k-1}$ and $G_k^{(k-1)}$ are both ordinary. Now define $G^{(k)}$
be the dominant transform of $G^{(k-1)}$ for $G\in\sG$. Define
$\sG^{(k)}=\{G^{(k)}:G\in\sG\}$; by Lemma~\ref{buildingsetlemma},
$\sG^{(k)}$ is ordinary and define $\sS^{(k)}$ to be the induced
arrangement of $\sG^{(k)}$. Since $\sG^{(k)}$ is ordinary, we see
that $\sS^{(k)}$ is ordinary. Finally for $k=N$ we get $X_N=X_\sG$.

We note that the theorem includes the compactification scheme constructed in
\cite{hu03} as a special case. The fact that this scheme arises from a suitable
building set is checked in \cite{li06}.
\end{proof}

\begin{proof}[Proof of \ref{cor:main}]
To deduce the Corollary~\ref{cor:main} from Theorem~\ref{the:main}
it suffices to produce ordinary building sets to construct
$X[n],X\langle n\rangle,X^\Gamma$ etc. The building sets for these
are constructed in \cite{li06}. These building sets are building
sets of $X^n$. To prove that they are ordinary building sets if $X$
is ordinary, we note that the building sets for (1)-(3) consists of
diagonals or polydiagonals, i.e. self-products of $X$ embedded in
$X^n$ by various diagonals. Thus the ordinarity of these building
sets follows from Proposition~\ref{ekedahl} by ordinarity of
self-products of ordinary varieties. For constructing $X^n_D$, we
use a building set which is constructed by \cite{kim08}, from $X^n$,
by blowing up a suitable subschemes which are self products of
$D,X$. By Proposition~\ref{ekedahl} this gives an ordinary building
set. The result follows from Theorem~\ref{the:main}. To construct
$X_D[n]$, we start with an ordinary building set in $X^n_D$,
consisting of the proper transform in $X^n_D$ of the multi-diagonals
in $X^n$. This is again an ordinary building set.

(5) This assertion is strictly part of the formalism of wonderful
compactification via \cite{kapranov93} (see \cite{li06}) but may be
of independent interest and so we give a proof for the sake of
completeness using \cite{keel92} where $\Mb_{0,n}$ is constructed
from $\P^1\times\cdots\times\P^1$ by a suitable sequence of blowups
with smooth, ordinary centers which are related to $\M_{0,j}$ for
$j<n$. In \cite{keel92} provided a construction of $\Mb_{0,n}$ as a
sequence of blowups and products. We will prove
Theorem~\ref{cor:main}(5) by induction on $n$.  Suppose $n=3$ then
$\Mb_{0,n}$ is a point hence is ordinary. Assume that $n=4$, then
$\Mb_{0,4}=\P^1$ hence is ordinary. Assume that the ordinarity of
$\Mb_{0,j}$ has been established for some all $j\leq n$; we will
show that $\Mb_{0,n}$ is also ordinary. Recall the construction of
\cite{keel92} (we will notations of that paper for this proof). We
let $B_1=\Mb_{0,n}\times\Mb_{0,4}=\Mb_{0,n}\times\P^1$. Then as
$\Mb_{0,n}$ is ordinary by induction and as $\P^1$ is ordinary, so
we deduce that $B_1$ is ordinary (see Proposition~\ref{ekedahl}).
The construction of $\Mb_{0,n}$ shows that for each subset
$T\subset\left\{1,2\ldots,n\right\}$ with $|T^C|\geq 2$, there
exists a collection of smooth subschemes $D^T$. For each $T$ these
subschemes are isomorphic to $\Mb_{0,i}\times\Mb_{0,j}$ for suitable
$i,j<n$. Thus by our induction hypothesis and
Proposition~\ref{ekedahl} we see that $D^T$ are ordinary and hence,
by Proposition~\ref{illusie} so is the blowup of $B_1$ along these
$D_T$ for every $T$. Thus $B_2$ is ordinary. More generally $B_k\to
B_{k-1}$ is the blowup of $B_{k-1}$ along the (disjoint) union of
strict transforms of  $D^T$ (for $|T^C|=k+1$) under $B_k\to B_1$.
Then $B_k$ is ordinary as $D^T$ are isomorphic to
$\Mb_{0,i}\times\Mb_{0,j}$ for suitable $i,j<n$. Thus $B_k$ is
ordinary and $\Mb_{0,n+1}=B_{n-2}$. This proves the assertion.

For (6) this is not immediate from \cite{li06} so we recall that
$T_{d,n}$ is constructed in \cite[Theorem~3.3.1]{chen06} in a manner
similar to the Fulton-MacPherson configuration scheme $X[n]$. The
procedure is inductive and starting from
$T_{1,3}=\P^1,T_{d,2}=\P^{d-1}$, $T_{1,n}=\Mb_{0,n+1}$ (note that by
the previous results these are all ordinary) we construct $T_{d,n}$
as follows: suppose $T_{d,n}$ has been constructed for some $d,n$.
Then $T_{d,n+1}$ is a sequence of blowups of a projective bundle
over $T_{d,n}$. Since the later is ordinary by induction, so is the
projective bundle over $T_{d,n}$ (by Proposition~\ref{illusie}). The
next blowups are along subschemes of the projective bundle which can
be identified with $T_{d,i}\times T_{d,j}$ for $j<n+1$ and so these
subschemes are ordinary by Proposition~\ref{ekedahl}. This proves
the assertion.
\end{proof}

\begin{remark}\label{biswas-rem}
Indranil Biswas has asked us the following question: if $X$ is a
smooth, projective ordinary surface, then is ${\rm Hilb}_n(X)$
ordinary for all $n\geq 1$? We note that it is known that if
$X$ is Frobenius split,  smooth, projective surface by then
\cite{kumar01}  ${\rm Hilb}_n(X)$ is Frobenius split.
By \cite{joshi03} smooth, proper Frobenius split
surfaces are ordinary. However by \cite{joshi03} the class of Frobenius
split varieties is not a subclass of ordinary varieties in
dimensions at least three and we note that the class of ordinary
surfaces is much bigger--for instance it includes general type
surfaces in $\P^3$ by the result of \cite{illusie90}. In any case Biswas'
question presents a natural variant of \cite{kumar01}.

Unfortunately we do not know how to answer Biswas's question.
The methods outlined here are not adequate as they require a far
better understanding of the geometry of ${\rm Hilb}_n(X)$ than we
seem to have at the moment. We note however that we can easily
deduce the result for ${\rm Hilb}_2(X)$ from our result for the
Fulton-MacPherson configuration space $X[2]$.
\end{remark}

\providecommand{\bysame}{\leavevmode\hbox to3em{\hrulefill}\thinspace}
\providecommand{\MR}{\relax\ifhmode\unskip\space\fi MR }
\providecommand{\MRhref}[2]{%
  \href{http://www.ams.org/mathscinet-getitem?mr=#1}{#2}
}
\providecommand{\href}[2]{#2}

\end{document}